\documentclass[12pt]{amsart}
\usepackage{color}
\definecolor{light}{gray}{0.9}

\newtheorem{theorem}{Theorem}

\theoremstyle{definition}
\newtheorem{remark}[theorem]{Remark}

\def\GL{\operatorname{GL}}
\def\rank{\operatorname{rank}}
\def\Rad{\operatorname{Rad}}
\def\QQ{{\mathbb Q}}
\def\ZZ{{\mathbb Z}}
\def\CC{{\mathbb C}}
\let\tensor=\otimes
\let\phi=\varphi
\def\Tau{{\textup{T}}}

\def\Box#1#2#3{\multiput(#1,#2)(1,0){2}{\line(0,1)1}
                           \multiput(#1,#2)(0,1){2}{\line(1,0)1}
                           \put(#1,#2){\makebox(1,1){$#3$}}}
\def\EBox#1#2{\Box#1#2{}}
\begin{document}

\title{Minors of symmetric and exterior powers}
\author{Winfried Bruns \and Wolmer V. 
Vasconcelos}
\address{Universit\"at Osnabr\"uck, FB Mathematik/Informatik,
49609 Osna\-br\"uck, Germany}
\email{winfried@mathematik.uni-osnabrueck.de}
\address{Department of Mathematics, Rutgers University,
Piscataway, NJ 08854-8019, USA} \email{vasconce@math.rutgers.edu}
\begin{abstract}
We describe some of the
determinantal ideals attached to symmetric, exterior and tensor
powers of a matrix. The methods employed use elements of Zariski's
theory of complete ideals and of representation theory.
\end{abstract}
\maketitle

Let $R$ be a commutative ring. The determinantal ideals attached to
matrices with entries in $R$ play ubiquitous roles in the study of the
syzygies of $R$--modules. In this note, we describe some the
determinantal ideals attached to symmetric, exterior and tensor
powers of a matrix. The methods employed use elements of Zariski's
theory of complete ideals and of representation theory, the results being
sharper for  rings containing  the rationals. 

Let $R$ be an integral domain (or a field) and $\phi:R^m\to R^n$
an $R$-linear map of rank $r$. It is an easy exercise to show that
the $d$-th symmetric power $S^d(\phi):S^d(R^m)\to S^d(R^n)$ has
rank $\binom{r+d-1}{d}$. Let $I_t(\phi)$ denote the ideal
generated by the minors (of a matrix representing $\phi$). Since
$$
\rank \phi=\max\{r: I_r(\phi)\neq 0\},
$$
one can immediately determine the radicals of the ideals
$I_t(S^d(\phi))$, namely
$$
\Rad I_t(S^d(\phi))=\Rad I_r(\phi)\quad \text{if}\quad
\binom{r+d-1}{d}\le t < \binom{r+d}{d}.
$$
Here we want to derive a more precise description of
$I_t(S^d(\phi))$, the ideals of minors of the exterior powers of
$\phi$, and the ideals of minors of a tensor product.

At least for the ideals associated with $\rank\phi=n$ one has a
very satisfactory result. To simplify notation we set
$I(\phi)=I_n(\phi)$.

Suppose that $m=n$. Then the ideals $I(\phi)$ and $I(S^d(\phi))$
are principal, generated by the determinants of square matrices,
and
\begin{equation}
\det(S^d(\phi))=\det(\phi)^s,\qquad
s=\binom{n+d-1}{d-1},\label{det}
\end{equation}
as follows immediately by transformation to a triangular matrix
(it is enough to consider a generic matrix over $\ZZ$, which is
contained in a field). For non-square matrices we can replace
$\det(\phi)$ by $I(\phi)$:

\begin{theorem}\label{maximal}
Let $R$ be a commutative ring, and $\phi:R^m\to R^n$ and
$\psi:R^p\to R^q$ be $R$-linear maps. Set $s=\binom{n+d-1}{d-1}$
and $e=\binom{n-1}{d-1}$.
\begin{enumerate}
\item[(1)] Then
$$
I(\phi\tensor\psi)\subset I(\phi)^qI(\psi)^n,\qquad
I(S^d(\phi))\subset I(\phi)^s,\qquad
I(\textstyle{\bigwedge^d}(\phi))\subset I(\phi)^e,
$$
with equality up to integral closure.
\item[(2)]  Suppose that $R$ contains the field of
rational numbers. Then equality holds in \textup{(1)}.
\end{enumerate}
\end{theorem}

\begin{proof}
For simplicity of notation we only treat the case of the symmetric
powers in detail. That of the exterior powers is completely
analogous. The tensor product requires slight modifications, which
we will indicate below.

(1) The formation of $I(X)$ and $I(S^d(\phi))$ commutes with ring
extensions for trivial reasons. Thus it is enough to prove the
inclusion for $R=\ZZ[X]=\ZZ[X_{ij}: i=1,\dots,n,\ j=1,\dots,r]$
and the linear map $\phi:R^m\to R^n$ given by the matrix
$X=(X_{ij})$. It is well-known (Trung \cite{Tr}; also see Bruns
and Vetter \cite[(9.18)]{BV}) that the ideals $I(X)^k$ are primary
with radical $I(X)$. It follows that $I(X)^k$ is integrally
closed, and therefore
$$
I(X)=\bigcap I(X)V
$$
where $V$ runs through the discrete valuation rings extending $R$.
To sum up: we may assume that $R$ is a discrete valuation ring.

Note that all ideals under consideration are invariant under base
change in $R^n$ and $R^m$. In fact, they are Fitting invariants of
$\phi$ and $S^d(\phi)$. By the elementary divisor theorem we can
therefore assume that $\phi$ is given by a matrix with non-zero
entries only in the diagonal positions $(i,i)$, $I=1,\dots,n$.
Then $S^d(\phi)$ is also given by a diagonal matrix, and it is an
easy exercise to show that indeed $I(S^d(\phi))=I(\phi)^s$ for
such matrices $\phi$.

As we have observed, equality holds for the ideals under
consideration if $R$ is a discrete valuation ring. This implies
equality up to integral closure.

(2) We have to prove the inclusion opposite to that in (1), and it
is enough to prove it for $R=\QQ[X]$.

Since $\rank \phi=n$, the linear map $S^d(\phi)$ has rank
$\binom{n+d-1}{d}$. Therefore the ideal under consideration is
non-zero, and its generators have the same degree as those of
$I(\phi)^s$.

The group $G=\GL_m(\QQ)\times \GL_r(\QQ)$ acts on $R$ via the
linear substitution sending each entry of $X$ to the corresponding
entry of
\begin{equation}
AXB^{-1},\qquad A\in \GL_m(\QQ),\ B\in \GL_r(\QQ).\label{GL}
\end{equation}
It is well-known that the $\QQ$-vector space $W$ generated by the
degree $rs$ elements of $I(X)^s$ is the irreducible
$G$-representation associated with Young bitableaux of rectangular
shape $s\times n$ (see De Concini, Eisenbud, and Procesi
\cite{DEP} or \cite[Section 11]{BV}). So the desired inclusion
follows if the ideal $I(S^d(\phi))$ is $G$-stable.

This is not difficult to see. In fact let, $\sigma$ be the
automorphism induced by the substitution \eqref{GL}. We write
$\sigma(\phi)$ for the linear map which we obtain from $\phi$ by
replacing each entry with its image under $\phi$, i.~e.\
$\sigma(\phi)=\phi\tensor \sigma$ for the ring extension
$\sigma:R\to R$. Then
\begin{align*}
\sigma\left(I(S^d(\phi))\right)&=I\left(\sigma(S^d(\phi))\right)=
I\left(S^d(\sigma(\phi))\right)\\
 &=I\left(S^d(\beta^{-1}\phi\alpha)\right)=
I\left(S^d(\beta)^{-1}\circ S^d(\phi)\circ
S^d(\alpha)\right)\\
 &=I(S^d(\phi))
\end{align*}
where $\alpha$ is the automorphism of $R^m$ induced by the matrix
$A$, and $\beta$ the automorphism of $R^n$ induced by $B$. The
very last equation uses again that Fitting ideals are invariant
under base change in the free modules.

For the tensor product the arguments above have to be modified at
two places. We can of course assume that $\psi$ is also given by a
matrix $Y$ of indeterminates. The first critical point is whether
$I(X)^qI(Y)^n$ is again integrally closed as an ideal of
$R=\ZZ[X,Y]$.

We can argue as follows: $\ZZ[X,Y]$ is a free $\ZZ$-module whose
basis is given by the products $\mu\nu$ where $\mu$ is a standard
monomial in the variables $X_{ij}$ and $\nu$ is a standard
monomial in the $Y_{kl}$ (see \cite[Section 4]{BV}). An element
$f$ of $\ZZ[X,Y]$ belongs to $I(X)^q$ if and only if $\mu\in
I(X)^q$ for all the factors $\mu$ appearing in the representation
of $f$ as a linear combination of standard monomials, and a
similar assertion holds with respect to $I(Y)^n$ and the factors
$\nu$. It follows immediately that $I(X)^qI(Y)^n=I(X)^q\cap
I(Y)^n$. Since both $I(X)^q$ and $I(Y)^n$ are integrally closed,
their intersection is integrally closed, too.

The other crucial question is whether the $\QQ$-vector space
generated by the degree $rq$ elements of $I(X)^qI(Y)^n$ is an
irreducible $G\times G'$-representation where $G'=\GL_p(\QQ)\times
\GL_q(\QQ)$ acts on $\QQ[Y]$ in the same way as $G$ on $\QQ[X]$,
and the action of $G\times G'$ on $K[X,Y]=K[X]\tensor K[Y]$ is the
induced one. The answer is positive since $W\tensor W'$ is an
irreducible $G\times G'$-representation if $W$ and $W'$ are
irreducible for $G$ and $G'$ respectively. It is enough to test
this after extending $\QQ$ by $\CC$, and then one can apply a
classical theorem (for example, see Huppert \cite[II.10.3]{Hup}).
\end{proof}

\begin{remark}\label{max-rem}
(a) One can formulate an abstract version of the theorem as
follows. Suppose that for each commutative ring $R$ one has a
functor $F_R$ on the category of $R$-modules satisfying the
following conditions:
\begin{itemize}
\item[(i)] $F_R$ commutes with ring extensions, i.~ e.\
$F_R(M\tensor S)=F_S(M\tensor R)$ for $R$-modules $M$ (and
similarly for $R$-linear maps) if $R\to S$ is a homomorphism of
rings;
\item[(ii)] $F_R(R^n)=R^{\lambda(n)}$ for all $n$;
\item[(iii)] if $V$ is a discrete valuation ring, then
$I(F_V(\phi))=I(\phi)^{\mu(m,n)}$ for a $V$-linear map $V^m\to
V^n$;
\end{itemize}
then both parts (1) and (2) Theorem \ref{maximal} hold
accordingly.

The Schur functors \cite{ABW}, which generalize symmetric and
exterior powers, satisfy these conditions.

(b) Especially the case of the symmetric powers is always
suspicious to depend on characteristic, but we have not been able
to find a counterexample to Theorem \ref{maximal}(2) in positive
characteristics.

(c) Suppose that $F_R(\phi)$ (as in (b)) is alternating when
$\phi$ (is given by an alternating matrix). Then it makes sense to
compare the Pfaffian ideals of $\phi$ (and $\psi$) with those of
$S^d(\phi)$, $\bigwedge^d(\phi)$, and $\phi\tensor\psi$.

For even $n$ (the matrix) $\phi$ has a Pfaffian, and one can
derive the formula analogous to equation \eqref{det} by
transformation to an anti-diagonal matrix. For odd $n$, and more
generally for lower order Pfaffians, one obtains variants of
Theorems \ref{maximal} and \ref{nonmax}, using arguments analogous
to those in the determinantal case. The neccessary representation
theory is contained in \cite{AF}.
\end{remark}

As we will see in the following, the situation is much more
complicated for lower order minors, and a description as precise
as Theorem \ref{maximal}(2) seems to be out of reach. Nevertheless
one can obtain reasonable upper and lower bounds from
specialization to diagonal matrices.

Let us recall some notation and facts from
\cite{BV} or \cite{DEP}. A \emph{Young diagram} (or partition) is a
non-increasing finite sequence $\sigma=(s_1,\dots,s_u)$ of
integers. We define the functions $\gamma_j$ on the set of all
Young diagrams by
$$
\gamma_j(\sigma)=\sum_{i=1}^u \max(0,s_i-j+1).
$$
The set of Young diagrams is partially ordered by
$$
\sigma\le\tau\quad\iff\quad \gamma_j(\sigma)\le
\gamma_j(\tau)\text{ for all } j.
$$
Furthermore one sets
$$
I^\sigma(\phi)=\prod_{i=1}^u I_{s_i}(\phi)\quad\text{and}\quad
I^{(\sigma)}(\phi)=\sum_{\tau\ge\sigma} I^\tau(\phi).
$$
If $X$ is a matrix of indeterminates (over some ring of
coefficients) we simply write $I^\sigma$ for $I^\sigma(X)$ etc. If
$\QQ\subset R$, then
\begin{equation}
I^{(\sigma)}=I^\sigma;\label{sigma}
\end{equation}
see \cite[(11.2)]{BV}. (This equation already holds if $m!$ or
$n!$ is a unit in $R$.)

Suppose $K$ is a field of characteristic $0$ and $X$ an $m\times
n$ matrix of indeterminates. Then $K[X]$ splits into a direct sum
$$
K[X]=\bigoplus_{\sigma} M_\sigma
$$
of pairwise non-isomorphic irreducible $G=\GL_m(K)\times
\GL_n(K)$-repre\-sentations ($\sigma$ runs through all Young
diagrams with at most $\min(m,n)$ columns). As a $G$-module, $M_\sigma$ is
generated by the doubly initial tab\-leaux of shape $\sigma$, i.~e.\
the product of the determinants of the $s_i\times s_i$-submatrices
of $X$ in its upper left corner, $i=1,\dots,u$. The decomposition
commutes with field extensions.

The third (and last) ideal associated with $\sigma$ is $I_\sigma$, 
the ideal generated by $M_\sigma$.
Evidently every $G$-stable ideal is a sum of ideals $I_\sigma$
(see \cite[Scetion 11]{BV}). If $\phi$ is a
matrix over a ring $R\supset\QQ$, then we may form
$I_\sigma(\phi)$ by substitution.

Suppose $f$ is a monomial in the indeterminates $Y_1,\dots,Y_n$.
Then $f$ has a unique decomposition $f=f_1\cdots f_u$ where each
$f_i$ is squarefree and $f_{i+1} \mid f_i$, $i=1,\dots,u-1$.
We set $s_i=\deg f_i$
and call $\sigma=(s_1,\dots,s_u)$ the \emph{shape} of $f$.

\begin{theorem}\label{nonmax}
Let $R$ be a ring, and $\phi:R^m\to R^n$ be an $R$-linear map.
Moreover, let $\Sigma$ be the set of the minimal elements (with
respect to $\le$) among the shapes of the monomials generating
$I_r(S^d(Y))$ for a diagonal matrix $Y$ of indeterminates (over
$\ZZ$).
\begin{enumerate}
\item[(1)] Then
$$
I_r(S^d(\phi))\subset \sum_{\sigma\in\Sigma} I^{(\sigma)}(\phi),
$$
with equality up to integral closure.
\item[(2)] If $\QQ\subset R$, then
$$
\sum_{\sigma\in\Sigma} I_\sigma(\phi)\subset I_r(S^d(\phi))\subset
\sum_{\sigma\in\Sigma} I^\sigma(\phi).
$$
\end{enumerate}
\end{theorem}

\begin{proof}
It is enough to prove the theorem for generic matrices $X$ with
$R=\ZZ[X]$ for (1) and $R=\QQ[X]$ for (2).

Let us first observe that the inclusion in (1) follows from the
same inclusion in (2), simply since $J=\sum_{\sigma\in\Sigma}
I^{(\sigma)}(X)$ has a basis of standard monomials, and $\ZZ[X]/J$
is therefore $\ZZ$-torsionfree.

Now we prove the inclusion $I_r(S^d(X))\subset J$ over $\QQ[X]$.
Note that $I_r(S^d(\phi))$ is $G$-stable, as discussed in the
proof of Theorem \ref{maximal}. Therefore it has a representation
$$
I_r(S^d(\phi))=\sum_{\tau\in\Tau} I_\tau.
$$
Evaluating the doubly initial tableau of shape $\tau$ on an
$m\times n$ diagonal matrix yields a monomial of shape $\tau$.
Thus $\tau\ge\sigma$ for some $\sigma\in\Sigma$ and, hence,
$I_\tau\subset I^{(\sigma)}$, as desired.

For the first inclusion in (2) we have to show that $\Sigma\subset
\Tau$. Observe that a monomial in $I^\tau(Y)$ always has shape
$\ge \tau$. Thus a monomial of shape $\sigma\in\Sigma$ can only be
contained in $I_r(S^d(Y))\subset \sum_{\tau\in\Tau} I^\tau(Y)$ if
there exists $\tau\in\Tau$ with $\tau\le\sigma$. But this implies
$\sigma\in\Tau$.

As in the proof of Theorem \ref{maximal} it is enough to prove the
assertion about integral closure for diagonal matrices of
indeterminates over $\ZZ$. But then we are comparing two monomial
ideals, and may pass to $\QQ$, whence it is enough to show
equality up to integral closure in (2). This however follows from
\cite[8.1]{DEP}: $I^\sigma$ is the integral closure of $I_\sigma$.
\end{proof}

In the theorem one can replace symmetric powers by exterior powers
and, more generally, establish a functorial version as indicated
in Remark \ref{max-rem}. A variant for tensor products can also be
given.

We demonstrate by an example that both inclusions in Theorem
\ref{nonmax}(2) can be strict simultaneously. Choose a $3\times
3$-matrix of indeterminates over $\QQ$ and set $d=3$ and $n=3$. By
inspection of the monomials that generate $I_3(S^3(Y))$ for a
diagonal matrix $Y$ of indeterminates one sees that, with the
notation of the figure, $\Sigma=\{\sigma_1,\sigma_2\}$ is the
right choice. We have computed the ideal $I_3(S^3(X))$ with
\textsc{Singular} \cite{Si}. Its component in the lowest degree
$9$ has dimension 5610. The figure shows all the degree 9 Young
diagrams $\ge \sigma_1$ or $\ge\sigma_2$. The dimensions of the
representations $M_\tau$ can be computed by the hook formula, and
it turns out that
$$
I_3(S^3(X))=I_{\sigma_1}+ I_{\sigma_2}+ I_{\sigma_5}+ I_{\sigma_9}
$$
since only the sum on the right hand side yields the correct
dimension in degree 9. (Instead of comparing dimensions one could
test which doubly initial Young tableau yield elements in
$I_3(S^3(X))$.)
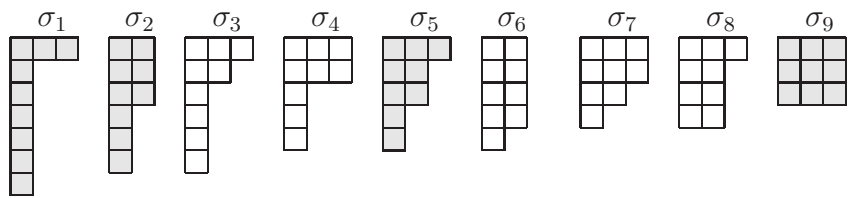
$$
\def\b#1#2{\EBox{#1}{#2}}\def\u{\unitlength}
\begin{picture}(3,8.5)(0,0)
{\color{light} \put(0,0){\rule{1\u}{7\u}}
\put(1,6){\rule{2\u}{1\u}} }
 \b00 \b01 \b02 \b03 \b04
 \b05 \b06 \b16 \b26 \put(1.2,7.5){$\sigma_1$}
\end{picture}\quad
\begin{picture}(2,6)(0,-1)
{\color{light} \put(0,0){\rule{1\u}{6\u}}
\put(1,3){\rule{1\u}{3\u}} }
 \b00 \b01 \b02 \b03 \b04 \b05 \b13
\b14 \b15 \put(0.7,6.5){$\sigma_2$}
\end{picture}\quad
\begin{picture}(3,6)(0,-1)
\b00 \b01 \b02 \b03 \b04 \b05 \b14 \b15 \b25
\put(1.2,6.5){$\sigma_3$}
\end{picture}\quad
\begin{picture}(3,6)(0,-2)
\b00 \b01 \b02 \b03 \b04 \b13 \b23 \b14 \b24
\put(1.2,5.5){$\sigma_4$}
\end{picture}\quad
\begin{picture}(3,6)(0,-2)
{\color{light} \put(0,0){\rule{1\u}{5\u}}
\put(1,2){\rule{1\u}{3\u}} \put(2,4){\rule{1\u}{1\u}}}
 \b00 \b01 \b02 \b03 \b04 \b12 \b13
\b14 \b24 \put(1.2,5.5){$\sigma_5$}
\end{picture}\quad
\begin{picture}(3,6)(0,-2)
\b00 \b01 \b02 \b03 \b04 \b11 \b12 \b13 \b14
\put(0.7,5.5){$\sigma_6$}
\end{picture}\quad
\begin{picture}(3,6)(0,-3)
\b00 \b01 \b02 \b03 \b11 \b12 \b13 \b22 \b23
\put(1.2,4.5){$\sigma_7$}
\end{picture}\quad
\begin{picture}(3,6)(0,-3)
\b00 \b01 \b02 \b03 \b10 \b11 \b12 \b13 \b23
\put(1.2,4.5){$\sigma_8$}
\end{picture}\quad
\begin{picture}(3,6)(0,-4)
{\color{light} \put(0,0){\rule{3\u}{3\u}} }
 \b00 \b01 \b02 \b10 \b11 \b12 \b20
\b21 \b22 \put(1.2,3.5){$\sigma_9$}
\end{picture}
$$


\begin{thebibliography}{1.}

\bibitem{AF} S. Abeasis and A. del Fra, \emph{Young diagrams and
ideals of Pfaffians}, Adv. Math. \textbf{35} (1980), 158--178.

\bibitem{ABW} K. Akin, D. A. Buchsbaum and J. Weyman, \emph{Schur functors
and Schur complexes}, Adv. Math. {\bf 44} (1982), 207--278.

\bibitem{BV}
W. Bruns and U. Vetter, \emph{Determinantal rings}, Springer
Lecture Notes {\bf 1327} (1988).

\bibitem{DEP} C. De Concini, D. Eisenbud, and C. Procesi, \emph{Young
diagrams and determinantal varieties}, Invent. math. \textbf{56}
(1980), 129--165.

\bibitem{Hup} B. Huppert, \emph{Endliche Gruppen I}, Springer
1967.

\bibitem{Si} G.-M. Greuel, G. Pfister, and H. Sch\"onemann. \textsc{Singular}
2.0. A Computer Algebra System for Polynomial Computations. Centre
for Computer Algebra, University of Kaiserslautern (2001). {\tt
http://www.singular.uni-kl.de}.


\bibitem{Tr} Ng\^o Vi\d{\^e}t Trung, \emph{On the symbolic powers of
determinantal ideals},  J. Algebra \textbf{58} (1979), 361--369.

\end{thebibliography}
\end{document}